\documentclass{amsart}\usepackage{amsmath,amssymb}
\usepackage[active]{srcltx} % ajout Kristian
\usepackage{mathrsfs}
\usepackage{pstricks}
\usepackage{pst-plot}
\usepackage[T1]{fontenc}
\usepackage{ifsym}
\usepackage{hyperref}

\usepackage{color}
\usepackage[colored]{shadethm}

\usepackage{graphics}
\usepackage[pdftex]{graphicx}

\usepackage{multicol}

\usepackage{pifont}

\definecolorseries{test}{rgb}{last}{yellow!50}{blue}

\let\RightarrowOrig\Rightarrow
\usepackage{marvosym}
\let\Rightarrow\RightarrowOrig

\newtheorem{definition}{\bf{Definition}}
\newtheorem{corollary}{\bf{Corollary}}
\newtheorem{theorem}{\bf{Theorem}}
\newtheorem*{theorem1}{\bf{Theorem 1}}
\newtheorem*{theorem2}{\bf{Theorem 2}}
\newtheorem*{theorem3}{\bf{Theorem 3}}

\newtheorem{lemma}{\bf{Lemma}}
\newtheorem{proposition}{\bf{Proposition}}
\newtheorem*{bohrtheo}{Classical Bohr's Theorem} 
\newtheorem*{Theorem}{Theorem}

\begin{document}
\LARGE
\title{  Estimates for the Bohr radius\\ of a Faber-Green condenser\\ in the complex plane.\\ \quad
\\ 
 \today}
\author{  {  P. Lass\`ere and E. Mazzilli}}
\date{\today}

\begin{abstract} We give some upper and lower estimates for the Bohr radius of a Faber-Green condenser in the complex plane. 
\end{abstract}


\address{LassËre Patrice: Institut de MathÈmatiques,  UMR CNRS 5580,
Universit\'e Paul Sabatier, 
118 route de Narbonne, 31062 TOULOUSE.}
\email{lassere@math.ups-tlse.fr}

\address{Emmanuel Mazzilli: UniversitÈ Lille 1, Villeneuve d'Ascq, 59655 Cedex.}
\email{Emmanuel.Mazzilli@math.univ-lille1.fr}  

\maketitle

\begin{quotation}
\begin{flushright}
\textit{``In memory of Professor Nguyen Thanh Van'',}
\end{flushright}
 
\end{quotation}
 
 \  \vspace{1cm}
 \vspace{1cm}
\section{Introduction}
 
\vspace{1cm}
The aim of this paper is to give some estimates for the Bohr radius of a Faber-Green condenser. Let us recall the classical Bohr Theorem for the unit disk: 

\bigskip
\begin{bohrtheo} \cite{bohr} Let $f(z)=\sum_n\,a_n z^n$ be holomorphic on the unit disc $\mathbb D:=\{\,z\in\mathbb C\ :\ \vert z\vert<1\,\}$. 
If $\vert f(z)\vert <1$ for all $z\in\mathbb D$, then $\sum_n\,\vert a_n\vert\cdot\vert z^n\vert<1$ for all $\vert z\vert<1/3$. 
Moreover, for all $\varepsilon>0$ there exists a holomorphic function  $f_\varepsilon(z)=\sum_n\,a^\varepsilon_n z^n$
  on $\mathbb D$ satisfying $\vert f_\varepsilon(z)\vert <1$ for all $z\in\mathbb D$, but $\sum_n\,\vert a^\varepsilon_n\vert\cdot\vert z^n\vert>1$ on $\vert z\vert= \varepsilon+1/3$. 
\end{bohrtheo}

\bigskip
For the last twenty years, this result has been generalized in many ways: to polynomials in one complex variable by Guadarrama \cite{GU}, Fournier \cite{F} and Chu \cite{CHU}. To several complex variables by Boas-Khavinson \cite{bk}, to the polydisk by 
Defant-Ortega-Cerd‡-OunaÔes-Seip, \cite{def}. To complex manifolds by Aytuna-Djakov \cite{ad} and by Aizenberg-Aytuna-Djakov in functional analysis \cite{AAD}. 
By Dixon \cite{dix}, Paulsen-Vern-Popescu-Singh \cite{paul} to   operator algebras. For a survey of literature on Bohr's phenomenon, see 
BÈnÈteau-Dahlner-Khavinson  \cite{BDK}.

\bigskip
In this paper, we focus on the Bohr radius of a condenser  in the complex plane.
For the convenience of the reader, let us recall the definition  introduced in \cite{lm1} (see Kaptanoglu-Sadik \cite{kap} for a partial approach in their seminal work).

 Denote $\mathbb D_r:=\{\,z\in\mathbb C\ :\ \vert z\vert <r\,\}$
and $\mathscr O(\mathbb D_r)$ the space of holomorphic functions on $\mathbb D_r$. We can reformulate the classical Bohr Theorem in the following way.

\bigskip
\textit{``The real $3$ is the smallest $r>1$ such that: if $f(z)=\sum_n\,a_nz^n\in\mathscr O(\mathbb D_r)$ , 
$\vert f(z)\vert<1$ on $\mathbb D_r$, then $\sum_n\,\vert a_n\vert\cdot\vert z^n\vert<1$ for all $z\in\mathbb D$.''}

\bigskip
This approach can be easily generalized for an arbitrary continuum (we recall that a continuum $K\subset\mathbb C$ is a compact set in $\mathbb C$ that contains at least two points and such that $\overline{\mathbb C}\setminus K$ is simply connected) if we notice that the discs  $\mathbb D_r$ are, for $r>1$, the levels sets of the Green function with pole at $\infty$ of $\overline{\mathbb C}\setminus \overline{\mathbb D}$. 

\bigskip
Given a continuum $K\subset\mathbb C$, by the Riemann mapping theorem, $\overline{\mathbb C}\setminus K$ has a Green function $\Phi_K$ 
with pole at $\infty$ and level sets $(\Omega^K_r)_{r>1}$. The sets $(K,(\Omega^K_r)_{r>1})$ will be called a Green-condenser. 
To achieve the construction, we have to ensure two things. First, we need to replace the Taylor 
basis $(z^n)_{n\geq 0}$ by a common basis $(\varphi_n)_{n\geq 0}$ for the spaces $\mathscr O(\Omega^K_r)$ 
(thanks to the general theory of common bases, there are many, \cite{ntv}) equipped with the usual compact convergence topology. 
We then consider a Green-condenser $(K,(\Omega^K_r)_{r>1}, (\varphi_n)_{n\geq 0})$ where $(\varphi_n)_{n\geq 0}$ is a common basis 
for the spaces $\mathscr O(\Omega^K_r)$. Second, we will use the following result (\cite{lm1} and \cite{ad}): 

\bigskip
\begin{Theorem}  
For a Green-condenser $(K,(\Omega^K_r)_{r>1}, (\varphi_n)_{n\geq 0})$, there always exists $r>1$ such that if $f=\sum_n\,a_n\varphi_n\in\mathscr O(\Omega_r^K)$ 
satisfies $\vert f\vert<1$ on $\Omega_r^K$, then 
$\sum_n\,\vert a_n\vert\cdot\Vert\varphi_n\Vert_K<1$.
\end{Theorem}

\bigskip
Note that in fact we obtained the result with the additional hypothesis that there exists $a\in K$ such that $\varphi_n(a)=0$ for all $n\geq 1$ and in \cite{ad}, 
Aytuna and Djakov relax this assumption even in a more general context. We can now define the Bohr radius for any condenser.

\bigskip
\textit{``The \textbf{Bohr radius} $B(K)$  of $(K,(\Omega^K_r)_{r>1}, (\varphi_n)_{n\geq 0})$ is the infimum of all $r>1$ such that $\Omega^K_r$ satisfies the previous theorem.''}

\bigskip
In the rest of the paper, we always work with $(F_{K,n})_{n\geq 0}$ the Faber basis for $K$ (see the definition in the next section)
and hence with the Faber-Green condenser $(K,(\Omega^K_r)_{r>1}, (F_{K,n})_{n\geq 0})$.  In general, 
it is not possible to calculate the exact value of $B(K)$ for an arbitrary continuum $K$. We know only the exact value of $B(K)$ 
in two cases: $K=\mathbb D$, of course, and for the elliptic condenser $K=[-1,1]$. Even in the elliptic case, the proof is difficult (see \cite{lm2}). 
The level sets $\Omega_r^{[-1,1]}$ of the Green function of $\overline{\mathbb C}\setminus [-1,1]$   are ellipses of loci $-1,1$ and 
eccentricity $\varepsilon=\frac{2r}{1+r^2}$ (the ``big level sets'' tend to ``big discs'' as $r\to\infty$). For this particular condenser, 
it is easy to deduce from \cite{lm2}  the exact value of $B(\Omega_r^{[-1,1]})$  for all $r>1$. Furthermore, we then can observe that $r\mapsto B(\Omega_r^{[-1,1]})$ is a decreasing function and tends to $3$ as $r$ tends to $\infty$. In fact, this last property is true for all condensers $(K,(\Omega^K_r)_{r>1}, (F_{K,n})_{n\geq 0})$ (see 
Theorem $2$). Let us point out that the classical fact ``big level sets'' tend to ``big discs'' as $r\to\infty$ 
is not enough to deduce this property. We have to analyze carefully, the behaviour of Faber polynomials and  of the Bohr radius $B(\overline{\Omega^K_r})$ 
for $r$ large (see the proof of Theorem 2).

\bigskip
In this paper, we give some estimates for  $B(K)$. The main results of the paper are:

\begin{theorem1} (uniform upper bound for $B(K)$. See section 3 for exact estimates)
\begin{enumerate}\item For every continuum $K\subset\mathbb C$, we have $B(K)\lesssim 13.8$.
\item Moreover, if $K$ is convex, then $B(K)\lesssim 5.26$.
\end{enumerate}
\end{theorem1}

\noindent\textbf{Remark: } If $K$ is the unit disk then $B(K)=3$, and if $K=[-1,1]$, $B(K)\simeq 5.1284$ (see \cite{lm2}).

\begin{theorem2} For every  continuum $K\subset\mathbb C$, we have
$$\lim_{r\to\infty}B(\overline{\Omega^K_r})=3.$$
\end{theorem2}
\
In a particular class of Faber-Green condenser (the positive class), we show the following result:

\begin{theorem3} For any positive Faber-Green condenser, we have
$$B(K)\geq 3.$$
Moreover, if $K$ is a positive Faber-Green condenser, then  $B(K)=3$ if and only if $K$ is a closed disk.

\end{theorem3}

The paper is organized as follows. 
The next section provides the background on Faber's polynomials. In section 3,  we prove Theorem 1 and some other estimates of $B(K)$ when $K$ is the interior of a  Jordan's curve or a $m$-cusped hypocycloid. Section 4 is devoted to the proof of the Theorem 2 and in the last section, we define the positive class for Faber-Green condenser and prove Theorem 3.

\vspace{1cm}
\section{Faber Polynomials}

\vspace{1cm}
This section is devoted to Faber polynomials and their properties. The classic reference on this topic is the book of P.K. Suetin \cite{suetin}.

\bigskip
First let us recall the construction of the Faber polynomials for a  continuum $K\subset\mathbb C$. Given a continuum $K\subset\mathbb C$, there exists a unique Riemann mapping
$$\Phi_K\ :\ \overline{\mathbb C}\setminus K\to 
\overline{\mathbb C}\setminus\mathbb D$$
normalized by 
$$\Phi_K(\infty)=\infty \quad\text{and}\quad\Phi'_K(\infty):=\lim_{z\to\infty}\frac{\phi_K(z)}{z}=\gamma>0,$$
where $\gamma$ is the logarithmic capacity or the transfinite diameter of $K$. In a neighborhood of the point $z=\infty$, we have the Laurent expansion:
$$\Phi_K(z)=\gamma z+\gamma_0+\frac{\gamma_1}{z}+
\frac{\gamma_2}{z^2}+\dots$$
Thus 
$$\begin{aligned}
\Phi_K^n(z)&=
\left( \gamma z+\gamma_0+\frac{\gamma_1}{z}+
\frac{\gamma_2}{z^2}+\dots\right)^n\\
&=\gamma^nz^n+a_{n-1}^{(n)}z^{n-1}+a_{n-2}^{(n)}z^{n-2}
+\dots+a_1^{(n)}z+a_0^{(n)}
+\sum_{j\geq 1}\,\frac{b_j^{(n)}}{z^j}.
\end{aligned}$$
The $n$-th Faber polynomial $F_{K,n}$ is now defined by taking the polynomial part of the Laurent expansion of $\Phi_K^n$.
For the sum of negative powers of $z$ we denote, 
$$E_{K,n}(z):=\sum_{j\geq 1}\,\frac{b_j^{(n)}}{z^j}=\Phi_K^n(z)-F_{K,n}(z).$$
If $R>1$ the circle $C(0,R)$ is mapped by $\Phi^{-1}_K$ onto a closed regular analytic curve $\Gamma_R$. This is the boundary of the bounded domain $\Omega^K_R=\{z\in\mathbb C\setminus K\ :\ \vert\Phi_K(z)\vert<R\,\}\cup K$ which is usually called the $R$-Green level set of $K$.

 \begin{center}
\includegraphics[scale=0.4]{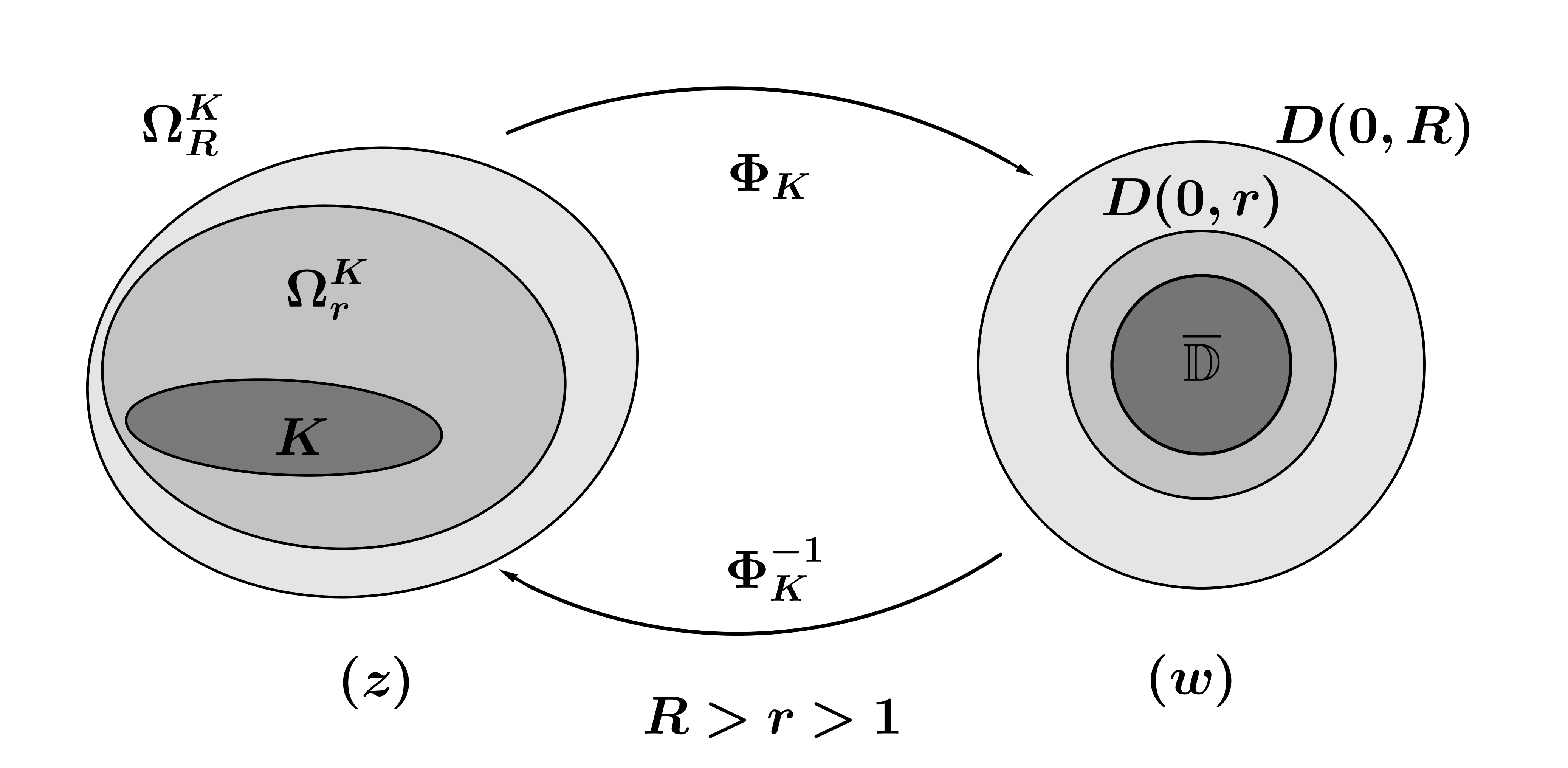}
\end{center}

Finally, the Faber-Green condenser is  $(K, (\Omega^K_R)_{R>1},(F_{k,n})_{n\geq 0})$. For a Faber-Green condenser, the situation is fairly like the Taylor one for the disk $(\overline{\mathbb D(0,1)}, (D(0,R))_{R>1},(z^n)_{n\geq 0})$ in the following way \textit{(\cite{suetin}, chapter 1)}: for all $f\in\mathscr O(\Omega_R^K)$, there exists a unique sequence $(a_n)_n$ of complex numbers  such that $f=\sum_n\,a_nF_{K,n}$ in $\mathscr O(\Omega_R^K)$ equipped with its natural compact convergence topology. Moreover, for $f\in\mathscr O(K)$ then $\limsup_n\vert a_n\vert^{1/n}=R^{-1}$ if and only if $R$ is the largest Green-level set such that $f\in\mathscr O(\Omega_R^K)$.

\bigskip 
\noindent Some examples: $\bullet$ If $K$ is the unit disk $\mathbb D$, then $\Phi_K(z)=z$. Hence, in this case the Faber polynomials coincide with the Taylor polynomials  $F_{K,n}(z)=z^n$ and the Faber-Green level sets are concentric disks $\Omega^K_R=D(0,R)$.\\
$\bullet$ For $K=[-1,1]$ we have $\Phi_K(z)=z+\sqrt{z^2-1},\ z\in\mathbb C\setminus K$ (where the branch of the square root is taken so that $\Phi_K'(\infty)=2$). In this example the Faber polynomials are the Chebyshev polynomials of the first kind $F_{K,n}=T_n$ and the level sets are ellipses.

\bigskip
We can also replace $K$ by one of its level sets $\overline{\Omega^K_R}$. It is not difficult to observe that
$$F_{K,n}(z)=R^nF_{\overline{\Omega^K_R},n}(z)$$
and we will often use this formula.

\bigskip
As usual when dealing with Faber polynomials,  it is better to work with the variable $w=\Phi_K(z)$ which lives in the annulus $D(0,R)\setminus \overline{\mathbb D}$ when $z\in\Omega^K_R\setminus K$. In this new coordinate one has:
 
 $$\begin{aligned} f(\Phi^{-1}_{K}(w))&=f(z),\quad \forall\,z=\Phi^{-1}_{K}(w)\in \Omega^K_R\setminus K\ \text{ and }\ w\in D(0,R)\setminus \overline{\mathbb D },\\
F_{K,n}(z)&=F_{K,n}(\Phi^{-1}_{K}(w))=w^n+\sum_{j\geq 1}^{\infty}{\alpha_{j}^{(n)}\over w^j}\\
\Phi^{-1}_K(w)&:={w\over \gamma}+\beta_0+\sum_{j\geq 1}^{\infty}{\beta_j\over w^j},\quad \forall\,\vert w\vert>1.
\end{aligned}$$

\vspace{1cm}
\section{Caratheodory-type inequalities and Uniform bounds}

\vspace{1cm}
\subsection{Caratheodory-type inequalities}

\bigskip
\begin{proposition} For all $R>1$ and $f=\sum_{n\geq 0}\,a_n F_{K,n}\in\mathscr O({\Omega^K_R})$ such that ${\texttt{re}}(f(z))\geq 0$, $z\in {\Omega^K_R}$, we have the Caratheodory-type inequalities:  \begin{equation}\vert a_n\vert\leq   \frac{ 2{\texttt{re}}(a_0)}{R^{n}-1},\quad\forall\,n\geq 1. \label{carateo}\end{equation}

\end{proposition}

 \bigskip
 \noindent\textbf{Proof: }  
$\bullet$ First, suppose that $f=\sum_{n\geq 0}\,a_n F_{K,n}\in\mathscr O(\overline{\Omega^K_R})$. We have, for all $1<r<R$,

\bigskip
\begin{equation}\int_{C_{(0,R)}}\,f(\Phi^{-1}_{K}(w))\,w^{n-1}dw
-\int_{C_{(0,r)}}f(\Phi^{-1}_{K}(w))w^{n-1}dw  =0.\label{inq1}\end{equation}

On the other way, because of the uniform convergence on compacts sets in $\overline{\Omega^K_R}\setminus K$:
$$\begin{aligned}
&\int_{C(0,R)}\,\overline{f(\Phi^{-1}_{K}(w))}\,w^{n-1}dw
-\int_{C(0,r)}\,\overline{f(\Phi^{-1}_{K}(w))}\,w^{n-1}dw\\
&=\sum_{j\geq 0}\left(\int_{C(0,R)}\,\overline{a_j}\cdot\overline{F_{K,j}(\Phi^{-1}_{K}(w))}\,w^{n-1}dw
-\int_{C(0,r)}\,\overline{a_j}\cdot \overline{F_{K,j}(\Phi^{-1}_{K}(w))}\,w^{n-1}dw\right).\\
\end{aligned}$$
But, $\displaystyle F_{K,j}(\Phi^{-1}_{K}(w))=w^j+\sum_{k\geq1}\frac{b_j^k}{w^k}$, thus 
$$\int_{C(0,r)}\,\overline{F_{K,j}(\Phi^{-1}_{K}(w))}\,w^{n-1}dw=\begin{cases}
0,&\quad\text{if }\ j\neq n,\\
2i\pi r^{2n}&\quad\text{if }\ j=n,\end{cases}$$
and
\begin{equation}\int_{_{C(0,R)}}\overline{f(\Phi^{-1}_{K}(w))}\,w^{n-1}dw
-\int_{C_{(0,r)}}\overline{f(\Phi^{-1}_{K}(w))}w^{n-1}dw
=2i\pi\overline{a_n}(R^{2n}-r^{2n}).\label{inq2}
\end{equation}
Then, \eqref{inq1}+\eqref{inq2} gives:
$$\begin{aligned}\int_{C(0,R)}\,2{\texttt{re}}({f(\Phi^{-1}_{K}(w))})w^{n-1}dw
&-\int_{C(0,r)}\,2{\texttt{re}}({f(\Phi^{-1}_{K}(w))})w^{n-1}dw\\
&=2i\pi \overline{a_n}(R^{2n}-r^{2n}).\end{aligned}$$
The real part of $f$ is positive for all $z\in\overline{\Omega^K_R}$, we get
$$\begin{aligned}\left\vert \overline{a_n}\cdot(R^{2n}-r^{2n})\right\vert&\leq \frac{1}{\pi}
\int_{C(0,R)}\,{\texttt{re}}({f(\Phi^{-1}_{K}(w))})\vert w^{n-1}\vert \cdot\vert dw\vert \\
&+\frac{1}{\pi}\int_{C(0,r)}\,{\texttt{re}}({f(\Phi^{-1}_{K}(w))})\vert w^{n-1}\vert \cdot\vert dw\vert\\
&=2{\texttt{re}}(a_0)(R^n+r^n)
\end{aligned}$$
and 
$$\vert a_n\vert\leq 2{\texttt{re}}(a_0)\frac{R^n+r^n}{R^{2n}-r^{2n}},\quad\forall\,n\geq 1,\ 1<r<R.$$
Let $r\to 1$:
$$\vert a_n\vert\leq 2{\texttt{re}}(a_0)\frac{R^n+1}{R^{2n}-1}= \frac{ 2{\texttt{re}}(a_0)}{R^{n}-1},\quad\forall\,n\geq 1.$$
$\bullet$ If $f\in\mathscr O(\Omega_R^K)$ then we get the previous inequality for every $R'<R$. It suffices then to take the limits when $R'$ goes to $R$.
\hfill$\blacksquare$

\bigskip
As a corollary we deduce the   estimates: 

\bigskip
\begin{theorem}  For every continuum  $K$, we have:
\begin{equation}  B(K)\leq \inf\left\{ R>1\ :\ \sum_{n\geq 1}\,\frac{4\sqrt{n\ln(n)+2n}}{R^n-1}\leq 1\right\}
\lesssim 13,8.\label{estim1}\end{equation}
  For every convex continuum $K$, we have:
\begin{equation} B(K)\leq \inf\left\{ R>1\ :\ \sum_{n\geq 1}\,\frac{4}{R^n-1}\leq 1\right\}\lesssim 5,26.\label{estim2}\end{equation}
 
\end{theorem}

 \bigskip
 \noindent\textbf{Proof: }Let $f=\sum_{n\geq 0}\,a_n F_{K,n}\in\mathscr O(\Omega_R^K)$ with $f(\Omega_R^K)\subset\mathbb D$. Up to a rotation, we can always suppose  $a_0\in\mathbb R^+$. Then the real part of $g:=1-f$ is positive on $\Omega_R^K$ and we can apply the Proposition 1 to 
 $$g(z)=1-a_0+\sum_{n\geq 1}a_n F_{K,n}(z).$$
This gives
$$\sum_{n\geq 0}\,\vert a_n\vert\cdot\Vert F_{K,n}\Vert_K\leq a_0+2(1-a_0)\sum_{n\geq 1}\,\frac{\Vert F_{K,n}\Vert_K}{R^n-1}.$$
So 
$$\sum_{n\geq 1}\,\frac{2\Vert F_{K,n}\Vert_K}{R^n-1}\leq 1 
\quad\implies\quad\sum_{n\geq 0}\,\vert a_n\vert\cdot\Vert F_{K,n}\Vert_K\leq 1$$
and $R\geq B(K)$. This gives immediately
$$B(K)\leq \inf\left\{ R>1\ :\ \sum_{n\geq 1}\,\frac{\Vert F_{K,n}\Vert_K}{R^n-1}\leq 1/2\right\}.$$
For $K$ convex we have: $1\leq \Vert F_{K,n}\Vert_K\leq 2$ (\cite{pom}). An easy computation gives
$$B(K)\leq \inf\left\{ R>1\ :\ \sum_{n\geq 1}\,\frac{4}{R^n-1}\leq 1\right\}\lesssim 5.26.$$
If $K$ is no more convex, then the sequence $(\Vert F_{K,n}\Vert_K)_n$ 
is no more bounded but cannot grow too fast. Crudely we have (\cite{pom})
$$1\leq \Vert F_{K,n}\Vert_K\leq 2\sqrt{n\ln(n)+2n}$$
and therefore
$$B(K)\leq  \inf\left\{ R>1\ :\ \sum_{n\geq 1}\,\frac{4\sqrt{n\ln(n)+2n}}{R^n-1}\leq 1\right\}.$$
Using Mapple, $B(K)\lesssim 13.8$. 
 \hfill$\blacksquare$
 
 \subsection{Example of the $m$-cusped hypocycloid}
\bigskip

The results of the  section 3 are particulary usefull when we can calculate exactly the norm of $F_{K,n}$. It is often the case when $$\Phi^{-1}_K(w):={w\over \gamma}+\beta_0+\sum_{j\geq 1}^{\infty}{\beta_j\over w^j}$$
where $\beta_j$, $j\geq 1$, are real and non-negative (in fact, it is the definition of the positive class of condenser see section 5). In this case, we have for the continuum $K$ (\cite{cur}, Theorem 3.1): 
$$F_{K,n}(\Phi^{-1}_K(w))=w^n+\sum_{j\geq 1}^{\infty}{\alpha_{j}^{(n)}\over w^j}\ \hbox{ where }\ \alpha_{j}^{(n)}\geq 0.$$
The $m$-cusped hypocycloids $H_m$ are in the positive class and satisfy the last property. Hypocycloids are starlike domains but not convex. Let us briefly recall the basic definitions and simple properties of the $m$-cusped hypocycloids. $H_m$ is the bounded region delimited by the closed curve $C_m$ defined by the equation
$$z=\exp(i\theta)+{1\over m-1}\exp(-(m-1)i\theta),\ \ m=2,3\cdots$$
The curve $C_m$ is the trajectory of a point on the unit disk  rolling without sliding in a larger disc of radius $m$. For $m=2$, $H_2=[-2,2]$ and we can calculate the exact value of $B([-2,2])$ (\cite{lm2}). If $m\geq 3$, it is straightforward to verify: $$\Phi^{-1}_{H_m}(w)=w+{1\over (m-1)w^{m-1}},$$ and $\Phi^{-1}_{H_m}$ admits a continuous extension on the unit circle which gives a topological mapping of the unit circle onto $C_m$. The coefficients $\alpha_{j}^{(n)}$ of $F_{H_{m,n}}$ are all positive and the series $\sum_{j}\alpha_{j}^{(n)}$ converge absolutely (see \cite{he}). This implies $$\Vert F_{H_m,n}\Vert_{H_m}= \vert F_{H_m,n}(\Phi^{-1}_{H_m}(1))\vert.$$
In \cite{he} it is implicitly proven that $$\Vert F_{H_3,n}\Vert_{H_3}=2+\left({-1\over 2}\right)^n:=M_{3,n}\text{\quad and\quad  }\Vert F_{H_4,n}\Vert_{H_4}=2+{\lambda^n+\bar\lambda^n\over 3^{n/ 2}}:=M_{4,n},$$ where $\lambda={1\over \sqrt{3}}(-1+\sqrt{2}i).$

\bigskip
We can now give the upper bound for $B(H_3)$ %and $B(H_4)$ 
using the same methods than in Theorem 1:

\bigskip
\begin{corollary} Let $i=3, 4$. Then for  $H_i$, we have the estimates:
$$B(H_i)\leq \inf\left\{R>1 : \sum_{n\geq 1}{2M_{i,n}\over R^n -1}\leq 1\right\}.$$
In particular : 
$$B(H_3)\leq 4.919167\dots$$
\end{corollary}

\bigskip
\noindent\textbf{Remark: } For $m>4$, we can prove the following
$$\Vert F_{H_m,n}\Vert_{H_m}\leq \left({m\over m-1}\right)^{m-1},$$ which is not optimal. The upper bound obtained for $B(H_m)$ with this estimate is not precise. 

\bigskip 

\subsection{Angular Measure}

\vspace{1cm}

In this subsection, we use classical notions (see \cite{pome} ). For the convenience  of the reader we recall these concepts here.
\bigskip

Suppose that $\Gamma$ is a rectifiable Jordan curve and let $\Omega$ be the interior of the bounded domain delimited by $\Gamma$. To almost every point of such curve $\Gamma$, we associate two angles as follows (see the picture below): 

\bigskip
\noindent$\bullet$ Let $s$ be the curvilinear coordinate of $\Gamma$. Then, for almost every $s$ we can define the tangent vector at $s$ to $\Gamma$. The first angle  $\sigma(s)$ will be the angle between the real axis and this tangent vector.

\bigskip

\noindent$\bullet$ For the second one, fix arbitrarily a point $z_0=\Phi(e^{i\varphi})\in \Gamma$. Define $v(\theta,\varphi)$ for $z=\Phi(e^{i\theta})\in \Gamma$ as $v(\theta,\varphi):=\texttt{arg}(z-z_0)$.

\begin{center}
\includegraphics[scale=0.11]{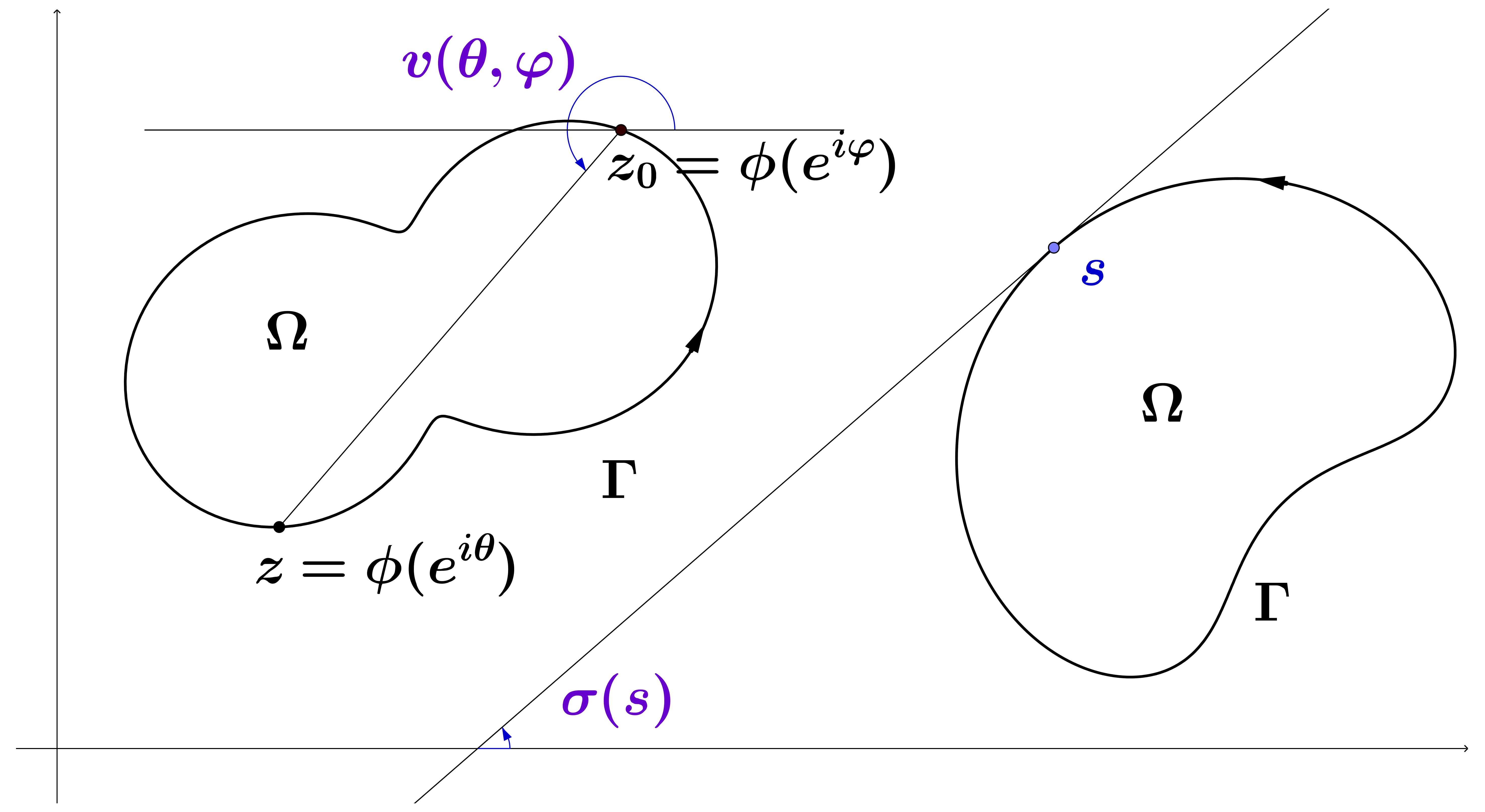}
\end{center}

If we suppose that $s\mapsto\sigma(s)$ is a function of bounded variation on $[0,l]$ ($l$ is the  length of $\Gamma$), we can associate to this function a unique finite total variation measure denoted also $\sigma$. Then we define 
$$V(\Gamma):=\int_0^l\,d\vert\sigma\vert (s).$$ 
If $\sigma$ is of bounded variation, then $\theta\mapsto v(\theta,\varphi)$ is also of bounded variation and it is not hard to see that 
$$\int_0^{2\pi}\,d\vert v\vert(\theta)\leq V(\Gamma).$$
\bigskip

We can then state the main result on the norm of Faber's polynomials using angle functions: 

\bigskip
\begin{proposition} (\cite{pome}). Suppose that $\Gamma$ is a rectifiable Jordan   curve and $\sigma$ of bounded variation then 
$$F_{\overline{\Omega},n}(\Phi(e^{i\varphi}))
=\frac{1}{\pi}\int_0^{2\pi}\,e^{in\theta}dv(\theta),$$
thus
$$ \Vert F_{\overline{\Omega},n} \Vert_{\overline{\Omega}}\leq \frac{V(\Gamma)}{\pi}.$$

\end{proposition}
We can give an other corollary of the Theorem 1: 

\begin{corollary}
Suppose $\Gamma$ and $\Omega$ are as before. Then, we have the estimate for the Bohr radius  
$$B(\overline{\Omega})\leq \inf
\left\{ \ R>1\ :\ \sum_{n\geq 1}\,\frac{2V(\Gamma) }{\pi(R^n-1)}\leq 1\ \right\}.$$
\end{corollary}
\bigskip
\noindent\textbf{Remarks: } The  quantity $V(\Gamma)$ is often easy to calculate or at least to estimate. Let us mention two examples: 
 \\ $\bullet$  If $\Gamma$ is convex, then obviously $V(\Gamma)=2\pi$. Note that, in this case, we get again the second part of Theorem 1.
\\ $\bullet$ If $\Gamma$ is a finite union of polygonal arcs, then the calculation of $V(\Gamma)$ is particulary simple. In the case of non convex such $\Gamma$, the angular approach gives a better estimate than  the general estimate of Theorem 1. 
 
% \vfill\eject
 \vspace{1cm}
 \section{Behaviour of $B(\overline{\Omega_r^K})$ when $r\to\infty$}

 \vspace{1cm}
  To simplify the notation when we consider $\overline{\Omega_{r}^{K}}$ as a continuum we will write $\Omega_{r}^{K}$. Thus, we write $F_{\Omega_{r}^{K},n}$ instead of $F_{\overline{\Omega_{r}^{K}}, n}$, the level sets $\Omega_{R}^{\Omega_{r}^{K}}$ instead of $\Omega_{R}^{\overline{\Omega_{r}^{K}}}$, $\Phi_{\Omega_{r}^{K}}$ instead of $\Phi_{\overline{\Omega_{r}^{K}}}$, and $B({\Omega_{r}^{K}})$ instead of $B(\overline{\Omega_{r}^{K}})$.
\bigskip

 Let $K\subset\mathbb C$ be a regular compact set, $r>1$. So $\Phi_{\Omega^K_r}=r^{-1}\Phi_K$ and therefore
 $$ F_{\Omega^K_r,n}(z)=\left(\frac{\Phi_K(z)}{r^n}\right)^n+\frac{E_n(z)}{r^n},\quad\forall\,z\in\mathbb C\setminus K,$$ 
or, in the $  w=r^{-1}\Phi_K(z) $ coordinate, 
\begin{equation} F_{\Omega^K_r,n}(\Phi_K^{-1}(rw))=w^n+\frac{E_n(\Phi_K^{-1}(rw))}{r^n},\quad\forall\,\vert w\vert>1.\label{tutu1}\end{equation}
Remember that (\cite{suetin}, pp.43), for all $1<r_0<r$, we have the uniform estimate:
\begin{equation}\vert E_n(z)\vert\leq \frac{r_0^n \texttt{length}(\partial\Omega^K_{r_0})}{2\pi \texttt{dist}(\partial\Omega^K_{r_0},\partial\Omega^K_r)},\quad \forall\,z\in \mathbb C\setminus\Omega^K_r,\,n\in\mathbb N,\label{tutu2}\end{equation}
where $\texttt{length}(\partial\Omega^K_{r_0})$ denotes the arclength of the level line $\{\vert \Phi_K\vert=r_0\}$.

\bigskip
Now let $0<r_1<1<r$ and $R>1$. Consider
$$f_{r_1}(z)=-r_1+\left(\frac{1}{r_1}-r_1\right)\sum_{n\geq 1}\frac{r_1^n}{R^n}F_{\Omega^K_r,n}(z).$$
Then $f_{r_1}\in\mathscr O(\Omega^{\Omega^{K}_r}_{R/r_1})\subset\mathscr O(\overline{\Omega^{\Omega^{K}_r}_R})$ and we have the following:

\begin{lemma} Let $r'_0>r_0>1$. There exists $M>0$  such that  
$$\sup_{z\in\partial\Omega^{\Omega^{K}_r}_R}\vert f_{r_1}(z)\vert\leq 1+M\left(\frac{1}{r_1}-r_1\right)\cdot\frac{1}{r},\quad \forall\,r>r'_0,\ 0<r_1<1.$$
\end{lemma}

\bigskip
\noindent\textbf{Proof: } The point $z\in\partial \Omega^{\Omega^{K}_r}_R$ if and only if $w=r^{-1}\Phi_{K}(z)=Re^{i\theta}$. So, using (\ref{tutu1}):
$$\begin{aligned}
f_{r_1}(\Phi_K^{-1}(rw))&=-r_1+\left(\frac{1}{r_1}-r_1\right)\sum_{n\geq 1}\frac{r_1^n}{R^n}\left(w^n+\frac{E_n(\Phi_K^{-1}(rw))}{r^n}\right)\\
&=-r_1+\left(\frac{1}{r_1}-r_1\right)\sum_{n\geq 1}\, r_1^n e^{in\theta} \\
&\qquad +\left(\frac{1}{r_1}-r_1\right)\sum_{n\geq 1}\,\frac{r_1^n}{R^nr^n} \cdot E_n(\Phi_K^{-1}(rw)) \\
&=\frac{e^{i\theta}-r_1}{1-r_1e^{i\theta}}+\left(\frac{1}{r_1}-r_1\right)\sum_{n\geq 1}\,\frac{r_1^n}{R^nr^n} \cdot E_n(\Phi_K^{-1}(rw))\\
&=\quad(A)\,\,+(B).
\end{aligned}$$
Since $r_1<1$
$$\Vert (A)\Vert=\sup_{\theta\in[0,2\pi]}
\left\vert\frac{e^{i\theta}-r_1}{1-r_1e^{i\theta}}\right\vert\leq 1.$$
For the second term,   (\ref{tutu2}) gives for all  $0<r_1<1<r_0<r_{0}^{'}<r$, $R>1$:
$$\begin{aligned}
\Vert(B)\Vert &\leq 
\left(\frac{1}{r_1}-r_1\right)\sum_{n\geq 1}\,\frac{r_1^n}{R^nr^n} \cdot \sup_{\vert w\vert=1}\vert E_n(\Phi_K^{-1}(rw))\vert\\
&\leq \left(\frac{1}{r_1}-r_1\right)\sum_{n\geq 1}\,\frac{r_1^nr_0^n}{R^nr^n} \cdot \frac{ \texttt{length}(\partial\Omega^K_{r_0})}{2\pi \texttt{dist}(\partial\Omega^K_{r_0},\partial\Omega^K_{r})}\\
&\leq \left(\frac{1}{r_1}-r_1\right)\cdot\frac{r_1r_0}{Rr}\cdot \frac{ \texttt{length}(\partial\Omega^K_{r_0})}{2\pi \texttt{dist}(\partial\Omega^K_{r_0},\partial\Omega^K_{r'_0})}\sum_{n\geq 0}\,\frac{r_1^nr_0^n}{R^nr^n}  \\
&\leq \frac{1}{r}\left(\frac{1}{r_1}-r_1\right)\cdot\frac{r_1r_0}{R}\cdot \frac{ \texttt{length}(\partial\Omega^K_{r_0})}{2\pi \texttt{dist}(\partial\Omega^K_{r_0},\partial\Omega^K_{r'_0})}\cdot\frac{1}{1-\frac{r_1r_0}{rR}} \\
& \leq \frac{1}{r}\left(\frac{1}{r_1}-r_1\right)\cdot\frac{r_1r_0}{R}\cdot \frac{ \texttt{length}(\partial\Omega^K_{r_0})}{2\pi \texttt{dist}(\partial\Omega^K_{r_0},\partial\Omega^K_{r'_0})}\cdot\frac{1}{1-\frac{r_0}{r'_0R}}\\
&\leq \frac{1}{r}\left(\frac{1}{r_1}-r_1\right)\cdot\frac{r_0r^{'}_0}{r^{'}_0-r_0}\cdot \frac{ \texttt{length}(\partial\Omega^K_{r_0})}{2\pi \texttt{dist}(\partial\Omega^K_{r_0},\partial\Omega^K_{r'_0})}\\
&= \frac{M}{r}\left(\frac{1}{r_1}-r_1\right),
\end{aligned}$$
where $M>0$ depends only of $r_0$ and $r^{'}_ 0$. 
Then 
$$\sup_{z\in\partial\Omega^{\Omega^{K}_r}_R}\vert f_{r_1}(z)\vert:=\Vert f_{r_1}\Vert_{\Omega^{\Omega^{K}_r}_R}\leq \Vert(A)\Vert+\Vert(B)\Vert\leq  1+M\left(\frac{1}{r_1}-r_1\right)\cdot\frac{1}{r},$$
for all $r,r_1,R$ such that $0<r_1<1,\ r>r'_0$ and $R>1$. 
\hfill$\blacksquare$

\bigskip
Suppose now $R>B(\Omega^K_r)$. The function  $f_{r_1}/\Vert f_{r_1}\Vert_{\Omega^{\Omega^{K}_r}_R}$ is holomorphic on  $\Omega^{\Omega^{K}_r}_R$ with values in $\mathbb D$, hence

$$r_1+\left(\frac{1}{r_1}-r_1\right)\sum_{n\geq 1}\frac{r_1^n}{R^n}\Vert F_{\Omega^K_r,n}\Vert_{\Omega^K_r}\leq \Vert f_{r_1}\Vert_{\Omega^{\Omega^{K}_r}_R}\leq 1+M\left(\frac{1}{r_1}-r_1\right)\cdot\frac{1}{r}$$
i.e.
$$\left(\frac{1}{r_1}-r_1\right)\sum_{n\geq 1}\frac{r_1^n}{R^n}\Vert F_{\Omega^K_r,n}\Vert_{\Omega^K_r}\leq
1-r_1+M\left(\frac{1}{r_1}-r_1\right)\cdot\frac{1}{r}.$$
Therefore
$$\sum_{n\geq 1}\frac{r_1^n}{R^n}\Vert F_{\Omega^K_r,n}\Vert_{\Omega^K_r}\leq
\frac{r_1}{1+r_1}+\frac{M}{r}.$$
With (\ref{tutu1}) and (\ref{tutu2}), we can write for any $r$ with $1<r_0<r$

$$\begin{aligned}\Vert F_{\Omega^K_r,n}\Vert_{\Omega^K_r}&\geq 1-\frac{\Vert E_n(\Phi_K^{-1}(rw))\Vert_{\partial\mathbb  D}}{r^n}\\
&\geq 1-\frac{r_0^n}{r^n}\cdot\frac{ \texttt{length}(\partial\Omega^K_{r_0})}{2\pi \texttt{dist}(\partial\Omega^K_{r_0},\partial\Omega^K_{r})}:=1-\frac{r_0^n}{r^n}\cdot M'(r).
\end{aligned}$$
In the same way, we have also an upper bound for $\Vert F_{\Omega_r,n}\Vert_{\Omega^K_r}$. So finally
\begin{equation}1-\frac{r_0^n}{r^n}\cdot M'(r) \leq \Vert F_{\Omega^K_r,n}\Vert_{\Omega^K_r}\leq 1+\frac{r_0^n}{r^n}\cdot M'(r).\label{norm}\end{equation}
Let $r$ such that  $ 1<r_0<r'_0<r$, then for any $R>B(\Omega^K_r)$,  we have
$$\begin{aligned}\sum_{n\geq 1}\frac{r_1^n}{R^n}\left(1-\frac{r_0^n}{r^n}\cdot M'(r) \right)&\leq \sum_{n\geq 1}\frac{r_1^n}{R^n}\Vert F_{\Omega^K_r,n}\Vert_{\Omega^K_r}\leq
\frac{r_1}{1+r_1}+\frac{M}{r},
\end{aligned}$$
for all  $0<r_1<1 $.
Thus
$$\begin{aligned}\frac{r_1}{R-r_1}&\leq \frac{r_1}{1+r_1}+\frac{M}{r}+M'(r)\cdot\sum_{n\geq 1}\,\left(\frac{r_0r_1}{Rr}\right)^n\\
&=\frac{r_1}{1+r_1}+\frac{M}{r}+M'(r)\cdot\frac{r_0r_1}{Rr-r_0r_1},\end{aligned}$$
for all $0<r_1<1$.
Letting now $r_1\to 1$, we get
$$\frac{1}{R-1}\leq \frac{1}{2}+\frac{M}{r}+M'(r)\cdot\frac{r_0}{Rr-r_0}
=\frac{1}{2}+\varepsilon(r),$$
where $\lim_{r\to\infty}\varepsilon(r)=0$ uniformly with respect to $R$ larger than one. \\
Hence for $R> B(\Omega^K_r)$ we have 
$$R\geq 3-\varepsilon'(r),$$
and so
$$B(\Omega^K_r)\geq 3-\varepsilon'(r),$$
where $\lim_{r\to\infty}\varepsilon'(r)=0$. Note that, in particular, 

\begin{equation}\liminf_{r\to+\infty} B(\Omega^K_r)\geq 3.\label{inf0} \end{equation}

\vspace{1cm}
\noindent $\bullet$ Now let us look for an upper bound for $B(\Omega^K_r)$ when $r$ is large. First observe that $( {\Omega^K_r}, (\Omega^K_{rR})_{R>1}, (F_{ {\Omega^K_r},n})_n)$ is the condenser associated with $ {\Omega^K_r}$ if $(K,(\Omega^K_r)_{r>1}, (F_{K,n})_n)$ is the condenser associated with $K$. 
\bigskip

Let $f=\sum_n\,a_n F_{  { \Omega^K_{r}},n}\in\mathscr O( {\Omega^{\Omega^K_r}_R})=\mathscr O(\Omega_{rR}^{K})$ be a function such that $f({\Omega^K_{rR}})\subset\mathbb D$. The proof of Proposition 1 on the annulus $A(\frac{1}{r},R)$ leads to 
$$\vert a_n\vert\leq \frac{2\texttt{re}(a_0)}{R^n-r^{-n}},\quad\forall\,n\in\mathbb N.$$
Assuming again $a_0\geq 0$, the Bohr phenomenon will occur if 
$$a_0+2(1-a_0)\sum_{n\geq 1}\,\frac{\left\Vert F_ { {\Omega^K_{r}},n}\right\Vert_{\Omega^K_r}}{R^n-r^{-n}}\leq 1.$$
This implies
$$2\sum_{n\geq 1}\,\frac{\left\Vert F_{  {\Omega^K_{r}},n}\right\Vert_{\Omega^K_r}}{R^n-r^{-n}}\leq 1.$$
From (\ref{norm}), it follows that  $\left\Vert F_{ {\Omega^K_r},n}\right\Vert_{\Omega^K_r}\leq 1+\frac{r_0^n}{r^n}\cdot M'(r)$ and thus
$$2\sum_{n\geq 1}\,\frac{\left\Vert F_{  {\Omega^K_{r}},n}\right\Vert_{\Omega^K_r}}{R^n-r^{-n}}\leq 2\sum_{n\geq 1}\,\frac{1+{r_0^n}{r^{-n}}\cdot M'(r)}{R^n-r^{-n}}.$$
If  $r>r_0>1$ is large enough, there exists a unique $R(r)>1$ such that
$$2\sum_{n\geq 1}\,\frac{1+ {r_0^n}r^{-n}\cdot M'(r)}{R(r)^n-r^{-n}}=1.$$
If $R_\infty:=\lim_{r\to\infty} R(r)$, we must have
$$2\sum_{n\geq 1} R_\infty^{-n}=1$$
that is $R_\infty=3$. On the other hand $R(r)\geq B(\Omega^K_{r})$ which implies, for $r$ large enough,  
\begin{equation}B(\Omega^K_{r})\leq 3+\varepsilon(r).\label{inf1}\end{equation}
Formulas (\ref{inf0}) and (\ref{inf1}) give 
$$ \lim_{r\to\infty}B(\Omega^K_{r})=3.$$
\hfill$\blacksquare$
%\vfill\eject

\quad

\bigskip

\bigskip
\section{The positive class of condenser and the proof of theorem 3}
\bigskip

Let us consider a special class of Faber-Green condenser:

\bigskip
\begin{definition}We say that $K$ is in the positive class of Faber-Green condenser or positive class, if we have for the continuum $K$:
$$z=\Phi^{-1}_{K}(w)={w\over \gamma} +\beta_0+\sum_{j=1}^{\infty}{\beta_j\over w^j},\ \hbox{with}\ \beta_j\geq 0, \forall\,j\geq 1.$$
\end{definition}

\bigskip
A continuum $K$ with this property has been considered by Curtiss and Pommerenke (\cite{cur}, \cite{pom}).
All the disks, all the lines, all the ellipses and all the $m$-cusped hypocycloids   are in this class. If $K$ is the closure of an analytic Jordan curve then $K$ being in the positive class implies that $K$ is a starlike domain (\cite{cur}, \cite{pom}). This class seems to be of some interest because we can evaluate precisely the sup-norm in $K$ of the Faber polynomials.

\bigskip
\noindent\textbf{Remark:} Clearly the Bohr radius is invariant by the automorphisms of the complex plane. Hence, the Theorem 3 is valid not only for the positive class but also for the pseudo-positive class: the orbit of the positive class by this group of automorphisms. Now the pseudo-positive class contains all the examples of continui considered by Eiermann and Varga in (\cite{ev}).

\subsection{Proof of the theorem 3}
\bigskip

Consider, for $r_1$ close to $1$, the family of functions
$$G_{r_1}={f_{r_1}\over \vert\vert f_{r_1}\vert\vert_{\Omega_{3}^{K}}}\ \hbox{with}\ f_{r_1}(z)=-r_1+\left({1\over r_1}-r_1\right)\sum_{n\geq 1}{r_{1}^n\over 3^n}F_{K,n}(z).$$
Clearly $G_{r_1}$  is holomorphic in $\overline{ \Omega_ {3}^{K}}$ and $\Vert G_{r_1}\Vert_{\Omega_{3}^{K}}\leq 1$. Suppose we have the Bohr property for $G_{r_1}$, i.e: $G_{r_1}:=\sum a_nF_{K,n}$ and $\sum\vert a_n\vert\cdot \Vert F_{K,n}\Vert_ {K}\leq 1$. This last inequality implies:
$$ r_1+\left({1\over r_1}-r_1\right)\sum_{n\geq 1}{r_{1}^{n}\over 3^n}\Vert F_{K,n}\Vert_K\leq \sup_{\vert w\vert =3}\left\vert f_{r_1}(\phi^{-1}_{K}(w))\right\vert\ \ (\bigstar).$$
$\bullet$ Estimate of $\vert f_{r_1}(\phi^{-1}_{K}(w))\vert$: 

In $w$-coordinate, $F_{K,n}(\phi^{-1}_K(w))=w^n+\sum_{j\geq 1}{\alpha_{j}^{(n)}\over w^j}$ with $\sum_{j\geq 1}{\alpha_{j}^{(n)}\over w^j}$ converges absolutely and uniformly on any compact set of $\{\vert w\vert >1\}$. For $w=3e^{i\theta}$, we have the following inequality
$$\begin{aligned}\vert f_{r_1}(\phi^{-1}_{K}(w))\vert & \leq \left\vert {e^{i\theta}-r_1\over 1-r_1e^{i\theta}}+\left({1\over r_1}-r_1\right)\sum_{n\geq 1}{r_{1}^n\over 3^n}\left (\sum_{j\geq 1}{\alpha_{j}^{(n)}\over 3^je^{i\theta j}}\right)\right\vert \\
 & \leq 1 + \left({1\over r_1}-r_1\right)\sum_{n\geq 1}{r_{1}^n\over 3^n}\left(\sum_{j\geq 1}{\alpha_{j}^{(n)}\over 3^j}\right).\end{aligned}$$
To get the previous inequality we use two facts: the double serie converges absolutely, and $\alpha_{j}^{(n)}$ are non-negative reals if $K$ is in the positive class. This is the crucial result of \cite{cur}.
\medskip

\noindent $\bullet$  Estimate of $\Vert F_{K,n}\Vert_K$:

Even in the positive class, $\sum_{j\geq 1}{\alpha_{j}^{(n)}\over w^j}$ is no longer abolutely convergent on $\{\vert w\vert =1\}$. We therefore have to modify the previous approach. Anyway, we have the following equality 
$$\Vert F_ {K,n}\Vert_{K}=\lim_{r\rightarrow 1^{+}}\sup_{\vert w\vert =r}\big\vert {F_ {K,n}(\phi^{-1}_ K(w))\over w^n}\big\vert.$$
 The term on the right is equal to $\lim_{r\rightarrow 1^{+}}(1+{1\over r^n}\sup_{\vert w\vert =r}\vert \sum_ {j\geq 1}{\alpha_{j}^{(n)}\over w^j}\vert)$ because $\alpha_{j}^{(n)}$ are non-negative reals by the Theorem of Curtiss\footnote{The hypothesis $K$ is in the positive class is crucial here to obtain a good lower bound for $\vert\vert F_ {K,n}\vert\vert_{K}$.}. Finally $$\vert\vert F_ {K,n}\vert\vert_{K}=1+\lim_{r\rightarrow 1^{+}}\sup_{\vert w\vert =r}\left\vert \sum_ {j\geq 1}{\alpha_{j}^{(n)}\over w^j}\right\vert\geq 1+\sum_{j\geq 1}{\alpha_{j}^{(n)}\over r^j},$$ for all $r>1$. Choose $1<r_0<3$, and suppose $(\bigstar)$ is valid. Then we must have
$$\left({1\over r_1}-r_1\right)\sum_{n\geq 1}{r_{1}^{n}\over 3^n}\left(1+\sum_{j\geq 1}{\alpha_{j}^{(n)}\over r_{0}^j}\right)\leq 1-r_1+\left({1\over r_1}-r_1\right)\sum_{n\geq 1}{r_{1}^{n}\over 3^n}\left(\sum_{j\geq 1}{\alpha_{j}^{(n)}\over 3^j}\right),$$
which implies the inequality
$${r_1\over 3-r_1}\leq {r_1\over 1+r_1}+\sum_{n\geq 1}{r_{1}^{n}\over 3^n}\left(\sum_{j\geq 1}\alpha_{j}^{(n)}\left({1\over 3^j}-{1\over r_{0}^j}\right)\right).$$
If $r_1$ tends to $1$, we obtain:
$$0\leq \sum_{n\geq 1}{1\over 3^n}\left(\sum_{j\geq 1}\alpha_{j}^{(n)}\left({1\over 3^j}-{1\over r_{0}^j}\right)\right).$$
Suppose $K$ is not a disk, then one of the $\alpha_{j}^{(n)}$ is strictly positive and the last inequality is not true. We have proved: if $K$ is not the disk then $B(K)>3$. We know by the classical Bohr Theorem that the Bohr radius of disks is $3$. The proof of Theorem $3$ is now complete.\hfill$\blacksquare$
\bigskip

As $H_m$ is in the positive class, we have the corollary:

\bigskip
\begin{corollary} If $H_m$ is the $m$-cusped hypocycloid then $B(H_m)>3$.
\end{corollary}

\bigskip
\section{Concluding remarks}
 
\bigskip

\begin{enumerate}\item  We are not able to prove the Theorem 3 in a larger class than the positive one. We suspect that the Theorem is true at least for starlike domains, but the proof seems delicate. Futhermore, it could be very interesting to produce a counter-example for general continuum.

\item For all convex continuum $K$, is it true or not that $B(K)\leq B([-1,1])$? Theorem 1 gives $B(K)\lesssim 5.26$ and in \cite{lm2} we proved that $B([-1,1])\simeq 5.1284$. The methods of Theorem 1 are far away from giving such inequality.
\bigskip

\item In general, it seems hopeless to compute the exact value of $B(K)$ for an arbitrary continuum $K$. But it seems possible to compute the exact value of the Bohr radius for the $3,4$-cusped hypocycloids $H_3$, $H_4$.

\end{enumerate}

\vfill
 
\end{document}